\documentclass[preprint,12pt,3p]{elsarticle}

\usepackage{amssymb}
\usepackage{amsmath,amsthm}
\usepackage{mathrsfs}
\usepackage{hyperref}
\usepackage{cleveref}
\usepackage[all,cmtip]{xy}
\usepackage{epsf, graphicx}
\usepackage{latexsym,amsfonts,amsbsy,amssymb}
\usepackage{geometry}
\usepackage{titletoc}
\usepackage{rotating}
\usepackage{algorithm}
\usepackage{algorithmic}
\usepackage{amscd}

\makeatletter

\@addtoreset{equation}{section} \makeatother
\newtheorem{theorem}{Theorem}[section]
\newtheorem{lemma}[theorem]{Lemma}
\newtheorem{remark}[theorem]{Remark}

\newtheorem{proposition}[theorem]{Proposition}

\def\O{\mathcal{O}}

\def\Ker{{\rm Ker}}
\def\Cok{{\rm Cok}}
\def\End{{\rm End}}

\def\I{{\rm Im}}
\def\F{\mathbb{F}}
\def\wt{\widetilde}

\def\SL{{\rm SL}(2,p^n)}
\def\GL{{\rm GL}(2,p^n)}

\makeatletter
\def\ps@pprintTitle{%
\let\@oddhead\@empty
\let\@evenhead\@empty
\def\@oddfoot{\reset@font\hfil\thepage\hfil}
\let\@evenfoot\@oddfoot
}
\makeatother

\begin{document}

\begin{frontmatter}

\title{The strengthened Brou\'{e} abelian defect group conjecture for ${\rm SL}(2,p^n)$ and ${\rm GL}(2,p^n)$}

\author[label1]{Xin Huang}
\address[label1]{SICM, Southern University of Science and Technology, Shenzhen 518055, China}
\address[label2]{Yau Mathematical Sciences Center, Tsinghua University, Beijing 100084, China}
\address[label3]{School of Mathematical Sciences, Peking University, Beijing 100871, China}


\author[label2]{Pengcheng Li}

\author[label1,label3]{Jiping Zhang}



\begin{abstract}
We show that each $p$-block of ${\rm SL}(2,p^n)$ and ${\rm GL}(2,p^n)$ over an arbitrary complete discrete valuation ring is splendidly Rickard equivalent to its Brauer correspondent, hence give new evidence for a refined version of Brou\'{e}'s abelian defect group conjecture proposed by Kessar and Linckelmann.

\end{abstract}


\end{frontmatter}


\section{Introduction}\label{s1}

Throughout this paper, $p$ is a prime number and $n$ is a positive integer. In representation theory of finite groups, the linear groups ${\rm SL}(2,p^n)$ and ${\rm GL}(2,p^n)$ are good examples for studying various theories. For instance, Bonnaf\'{e} \cite{Bonnafe} used the group $\SL$ to give an introduction to Harish-Chandra and Deligne-Lusztig theories. As noted in \cite{Bonnafe}, the group $\SL$ is sufficiently simple to allow a complete description, and yet sufficiently rich to illustrate some of the most delicate aspects of the theory.
Many global-local conjectures was verified for some blocks of the groups $\SL$ and $\GL$, such as Turull's refinement of the Alperin-McKay conjecture (see \cite[Theorem 4.11]{T08}) and Brou\'{e}'s abelian defect group conjecture (see \cite{Chuang},\cite{Oku},\cite{Yo},\cite{Marcus03}).

In this paper we investigate a refined abelian defect group conjecture for all $p$-blocks of $\SL$ and $\GL$. This refined conjecture was proposed by Kessar and Linckelmann (see \cite[page 186]{Kessar_Linckelmann}) and claimed that for any complete discrete valuation ring $\O$ and any block of a finite group over $\O$ with an abelian defect group, there is a splendid Rickard equivalence between the block algebra and its Brauer correspondent. Brou\'{e}'s original conjecture is with the assumption that the complete discrete valuation rings have splitting residue fields.

Throughout this paper, $k\subseteq k'$ are fields of characteristic $p$, $\O$ is a complete discrete valuation ring with residue field $k$, or $\O=k$. Assume that $k'$ is a splitting field for all finite groups considered below.

If $G$ is a finite group, we denote by $G^{\rm op}$ the opposite group, and we denote by $\Delta G$ the subgroup $\{(g,g^{-1})~|~g\in G\}$ of $G\times G^{\rm op}$.  By a {\it block} of the group algebra $\O G$, we mean a primitive idempotent $b$ of the center of $\O G$, and $\O Gb$ is called a {\it block algebra}. For a subgroup $H$ of $G$, let $(\O G)^H$ denote the set of $H$-fixed elements of the block algebra $\O G$ under the conjugation action. If $H$ is a $p$-subgroup, the {\it Brauer map} is the $\O$-algebra homomorphism
${\rm Br}_H: (\O G)^H\to kC_G(H)$, $\sum_{g\in G}\alpha_gg\mapsto \sum_{g\in C_G(H)}\bar{\alpha}_gg,$
where $\bar{\alpha}_g$ denotes the image of $\alpha_g$ in $k$.
For a block $b$ of $\O G$, a {\it defect group} of $b$ is a maximal $p$-subgroup $P$ of $G$ such that ${\rm Br}_P(b)\neq 0$. By Brauer's first main theorem, there is a unique block $c$ of $\O N_G(P)$ with defect group $P$ such that ${\rm Br}_P(b) = {\rm Br}_P(c)$ and the map $b\mapsto c$ is a bijection between the set of blocks of $\O G$ with defect group $P$ and the set of blocks of $\O N_G(P)$ with defect group $P$. This bijection is known as the {\it Brauer correspondence}.

Let $A$ and $B$ be symmetric $\O$-algebras. Let $X$ be a bounded complex of finitely generated $(A,B)$-bimodules which are projective as left $A$-modules and as right $B$-modules, and let $X^*:={\rm Hom}_\O(X,\O)$ be the dual complex. It is said that $X$ induces a {\it Rickard equivalence} and that $X$
is a {\it Rickard complex} if there exist a contractible complex of $(A,A)$-bimodules $Y$
and a contractible complex of $(B,B)$-bimodules $Z$ such that $X\otimes_B X^*=A\oplus Y$ as complexes of $(A,A)$-bimodules and $X^*\otimes_A X= B\oplus Z$ as complexes of $(B,B)$-bimodules.
 Let $G$ be a finite group and $H$ a subgroup of $G$. Let $b$ (resp. $c$) be an idempotent in the center of $\O G$ (resp. $\O H$). Let $X:=(X_n)_{n\in \mathbb{Z}}$ be a Rickard complex of $(\O Gb,\O Hc)$-bimodules. Recall that if every indecomposable direct summand of each $X_n$ is isomorphic to a direct summand of ${\rm Ind}_{\Delta H}^{G\times H^{\rm op}}(\O)$, then $X$ is said to be {\it splendid}; $\O Gb$ and $\O Hc$ are said to be {\it splendidly Rickard equivalent}.

The main result of this paper is the following.

 \begin{theorem}\label{main}
Let $G$ be $\SL$ or $\GL$, $b$ a block of $\O G$, then $\O Gb$ is splendidly Rickard equivalent to its Brauer correspondent algebra.
\end{theorem}

By lifting theorem of splendid Rickard equivalences (see Theorem \ref{lifting} below), to prove Theorem \ref{main}, we may assume that $\O=k$. We have seen in the proof of \cite[Theorem 1.12]{Kessar_Linckelmann} that in order to prove Theorem \ref{main}, it is enough to prove the following statement.

\begin{theorem}\label{main2}
Let $G$ be a $\SL$ or $\GL$, $b$ a block of $k'G$ having a defect group $P$. Let $c$ be the block of $k'N_G(P)$ corresponding to $b$ via the Brauer correspondence. Suppose that $b\in kG$. Then $c\in kN_G(P)$ and the block algebras $kGb$ and $kN_G(P)c$ are splendidly Rickard equivalent.
\end{theorem}

We prove this theorem for $\SL$ in Section \ref{proof} and for $\GL$ in Section \ref{proofforGL}.
It is well-known that if $p=2$, the group algebra $k'\SL$ has two blocks, i.e., the principal block and a defect zero block; if $p>2$, $k'\SL$ has three blocks, namely the principal block, a full defect non-principal block and a defect zero block. Okuyama \cite{Oku} proved that, if $k'$ is algebraically closed, the principal block algebra of $k'\SL$ is Rickard equivalent to its Brauer correspondent. Yoshii \cite{Yo} generalised Okuyama's method and proved the same result for the full defect non-principal block of $k'\SL$. Since $\F_{p^n}$ is a splitting field for both $\SL$ and the normaliser of a Sylow $p$-subgroup, all modules appeared in Okuyama and Yoshii's proofs are defined over $\F_{p^n}$. So Okuyama and Yoshii's proofs still work for replacing $k'$ by $\F_{p^n}$. We note that when Okuyama and Yoshii proved the existence of Rickard equivalences, they didn't point out whether the Rickard equivalences are splendid.

In \cite{Marcus03}, by establishing a graded version of Rickard's fundamental theorem, Marcus proved that, over a splitting field, there is a Rickard equivalence between the principal block algebra of $\GL$  and its Brauer correspondent algebra (see \cite[Example 3.14]{Marcus03}). The proof of Theorem \ref{main2} is based on Okuyama, Yoshii and Marcus' work.

\section{Block idempotents and coefficient rings}

The following theorem on lifting splendid Rickard equivalences is due to Rickard.

\begin{theorem}[{\cite[Theorem 5.2]{Rickardsplendid}}]\label{lifting}
Let $G$ be a finite group and $H$ a subgroup of $G$. Let $b$ (resp. $c$) be an idempotent in the center of $\O G$ (resp. $\O H$). Denote by $\bar{b}$ (resp. $\bar{c}$) the image of $b$ (resp. $c$) in $kG$ (resp. $kH$). Assume that there is a complex $\bar{X}$ of $(kG\bar{b}, kH\bar{c})$-bimodules inducing a splendid Rickard equivalence. Then there is a complex $X$ of $(\O Gb, \O Hc)$-bimodules inducing a splendid Rickard equivalence and satisfying $k\otimes_\O X\cong \bar{X}$.
\end{theorem}

Note that although the statement in \cite[Theorem 5.2]{Rickardsplendid} is for principal blocks, but
the proof carries over nearly verbatim to arbitrary blocks. We also note that the blanket
assumption in [12] that the coefficient rings are big enough is not used in the proof of \cite[Theorem 5.2]{Rickardsplendid}.

For notational convenience, we identify $\F_p$ and the prime field of any field of characteristic $p$. Let $G$ be a finite group. For $\alpha=\sum_{g\in G}\alpha_gg$ an element of $k'G$, denote by $k[\alpha]$ the smallest subfield of $k'$ containing $k$ and all coefficients $\alpha_g$, $g\in G$. If $\sigma$ is an automorphism of $k'$, then $\sigma$ induces a ring automorphism of $k'G$ (abusively still denoted by $\sigma$) in an obvious way. So $\sigma(\alpha)$ is $\sum_{g\in G}\sigma(\alpha_g)g$. The following proposition is well-known.

\begin{proposition}\label{2.2}
Let $G$ be a finite group and $b$ a block of $k'G$. Let $P$ be a defect group of $b$, and $c$ the block of $k'N_G(P)$ corresponding to $b$ via the Brauer correspondence.
Then $k[b]=k[c]$.
\end{proposition}

\noindent{\it Proof.} Since every finite group has a finite splitting field, we may assume that $k'$ is finite. Since $c\in k'C_G(P)$ (see e.g. \cite[Theorem 6.7.6 (ii)]{Linckelmann}), we have $c={\rm Br}_P(c)={\rm Br}_P(b)$, and hence $c\in k[b]$. For any $\sigma\in {\rm Gal}(k'/k[c])$, $\sigma(b)$ is also a block of $k'G$ with defect group $P$, and we have ${\rm Br}_P(\sigma(b))=\sigma({\rm Br}_P(b))={\rm Br}_P(c)$. By definition, $\sigma(b)$ is the Brauer correspondent of $c$, which yields $\sigma(b)=b$. So $b\in k[c]$. $\hfill\square$

\medskip Let $G$ be a finite group, $b$ a block of $k'G$, the smallest field $k$ such that $b\in kG$ is $\F_p[b]$. By the following proposition, to prove Theorem \ref{main2} for $\SL$, we may take $k=\F_p$.

\begin{proposition}\label{blockidempotentofSL}
Let $b$ be a block of $k'\SL$, then $\F_p[b]=\F_p$.
\end{proposition}

\noindent{\it Proof.} Since $\F_{p^n}$ is a splitting field of $\SL$, we may take $k':=\F_{p^n}$. The group ${\rm Gal}(k'/\F_p)$ acts via $\F_p$-algebra automorphisms on the group algebra $k'\SL$ in an obvious way. These automorphisms permute blocks, preserve defect groups of blocks, and fix the principal block. Recall that if $p=2$, the group algebra $k'\SL$ has two blocks, the principal block and a defect zero block; if $p>2$, $k'\SL$ has three blocks, the principal block, a full defect non-principal block and a defect zero block. So each block of $\F_{p^n}\SL$ is stable under the action of ${\rm Gal}(\F_{p^n}/\F_p)$, and hence contained in $\F_p\SL$. $\hfill\square$

\begin{lemma}\label{blockidempotentofGL}
Let $b$ be a block of $k'\SL$, then $b$ is an idempotent in the center of $k'\GL$. Let $\tilde{b}$ be a block of $k'\GL$ with $b\tilde{b}\neq 0$, then $b$ and $\tilde{b}$ have the same defect groups.
\end{lemma}

\noindent{\it Proof.} Since the elements of $\GL$ induce automorphisms of the $k'$-algebra $k'\SL$ by conjugation, the elements of $\GL$ permute the blocks of $k'\SL$.  These permutations preserve defect groups of blocks and fix the principal block. Hence each block of $k'\SL$ is stable under the action of $\GL$, the first statement holds. The second statement follows from \cite[\S15, Theorem 1 (2)]{Alperin}. $\hfill\square$

\section{Descent of tilting complexes arising from stable equivalences}\label{Section2}

Tilting complexes and split-endomorphism two-sided tilting complexes (which are also called Rickard complexes) were introduced by Rickard (\cite{RickardMorita}, \cite{Rickardderivedfunctor}, \cite{Rickardsplendid}). In this section, we review a way to construct a tilting complex via a stable equivalence of Morita type, discovered by Rouquier \cite{Rou} and improved by Okuyama \cite{Oku2}. These results were also reviewed in \cite[3.1]{Wong}. After that, we give a descent criterion, which is the main tool in proving Theorem \ref{main2}. Unless specified otherwise, modules in the paper are left modules. For an algebra $A$, $K^b(A)$ (resp. $D^b(A)$) denotes the homotopy (resp. derived) category of bounded complexes of finite generated $A$-modules; $A^{\rm op}$ denotes the opposite algebra.

\subsection{Tilting complexes}\label{tilting}
Let $A'$ and $B'$ be finite-dimensional symmetric $k'$-algebras, where $k'$ is a splitting field of both $A'$ and $B'$. Assume that a finitely generated $(A',B')$-bimodule $M'$ induces a stable equivalence of Morita type between $A'$ and $B'$. Assume that $A'$ and $B'$ have no semisimple summand and $M'$ has no projective summand. Let $\{T_i~|~i\in I\}$ be a set of representatives of classes of  simple $B'$-modules and $Q'_i$ be a projective cover of $T_i$. Let $P'_i$ be a projective cover of the $A'$-module $M'\otimes_{B'} T_i$. Then by \cite[Lemma 2]{Rou}, we can choose certain $(A',B')$-bimodule homomorphisms $\delta_i':P'_i\otimes_{k'} {Q'}_i^*\to M'$, such that $\oplus_{i\in I}\delta'_i:\oplus_{i\in I}P'_i\otimes_{k'} {Q'}_i^*\to M'$ is a projective cover of the $(A', B')$-bimodule $M'$.

Let $I_0$ be a fixed subset of $I$. Define the $(A',B')$-bimodule $P'(I_0)$ to be $\oplus_{i\in I_0}P'_i\otimes_{k'} {Q'}_i^*$ and denote by $\delta':=\delta'(I_0):P'(I_0)\to M'$ the restriction of $\oplus_{i\in I}\delta'_i$ to $P'(I_0)$. Define the complex $M'(I_0)^{\bullet}$ of $(A',B')$-bimodules to be
$$\cdots\to 0\to P'(I_0)\xrightarrow{\delta'}M'\to 0\to \cdots~,$$
where $M'$ lies in degree zero. By \cite[Theorem 1.1]{Oku}, $M'(I_0)^\bullet$ is a tilting complex for $A'$ if and only if a certain condition is satisfied.

\begin{theorem}\label{Ok1.2}
Keep the notation above. Assume that as a complex of projective $A'$-modules, $M'(I_0)^\bullet$ is a tilting complex for $A'$. Let $C':={\rm End}_{K^b(A')}(M'(I_0)^\bullet)^{\rm op}$. Then the $k'$-algebra $C'$ has $(B',B')$-bimodule structure induced from the right action of $B'$ on $M'(I_0)^\bullet$. The following holds.

\medskip\noindent{\rm(1) (\cite{RickardMorita}, \cite[Corollary 5.3]{Rickardderivedfunctor})}. $C'$ is a finite-dimensional symmetric $k'$-algebra, and is derived equivalent to $A'$.

\noindent {\rm (2) (\cite[page 5]{Oku}, see also \cite[Proposition 2.2.1]{Yo})}. The right action of $B'$ on $M'(I_0)^\bullet$ induces a $k'$-algebra monomorphism $\rho'$ from $B'$ to $C'$.

\noindent {\rm (3) (\cite[Theorem 1.2]{Oku})}. There exists a complex $N'(I_0)^{\bullet}$ of $(A',C')$-bimodules, which is a direct summand of $M'(I_0)^\bullet\otimes_{B'} C'$ such that $N'(I_0)^{\bullet}$ is a Rickard complex for $A'$ and $C'$.

\noindent {\rm (4) (\cite[Theorem 1.2]{Oku})}. $N'(I_0)^\bullet$ is of the form
$$\cdots\to0\to Q'\xrightarrow{f'} N'\to 0\to\cdots,$$
where $N'$ lies in degree $0$, $Q'$ is a projective $(A',C')$-bimodule.

\end{theorem}

\subsection{The construction of $f'$}\label{f'}
Now we begin to review \cite[(1.7)]{Oku} about the construction of the modules $Q'$, $N'$ and the homomorphism $f'$ in Theorem \ref{Ok1.2} (4). Let $Q'$ be an injective hull of the $(A',C')$-bimodule ${\rm Ker}\delta'$. Then by \cite[(1.7)]{Oku}, the sequence
$$0\to \Ker\delta'\to Q'\xrightarrow{\phi'} Q'/\Ker\delta'\to 0$$
is isomorphic to a direct summand of the sequence
$$0\to  \Ker\delta'\otimes_{B'} C'\to  P'(I_0)\otimes_{B'}C'\to \I\delta'\otimes_{B'}C'\to 0.$$
So we can write $\I\delta'\otimes_{B'}C'\cong Q'/\Ker\delta'\oplus Y'$ for some direct summand $Y'$ of $\I\delta'\otimes_{B'}C'$

Let $Q''$ be an injective hull of $Y'$. By \cite[(1.7)]{Oku}, the sequence
$$0\to Y' \to Q''\to Q''/Y'\to 0 $$
is isomorphic to a direct summand of the sequence
$$0\to \I\delta'\otimes_{B'}C'\to  M'\otimes_{B'}C' \to  \Cok\delta'\otimes_{B'}C'\to 0. \eqno(*)$$
Write $ M'\otimes_{B'}C'\cong Q''\oplus N''$ and $\Cok\delta'\otimes_{B'}C'\cong Q''/Y'\oplus W'$ for some appropriate $(A',C')$-bimodules $N''$ and $W'$. Then we have the resulting exact sequence
$$0\to Q'/\Ker\delta'\xrightarrow{\psi''}N''\to W'\to 0 \eqno(**)$$
as a direct summand of $(*)$.
By \cite[(1.7)]{Oku}, $\Cok \delta'$ is isomorphic to a direct summand of $W'$; when writing $W'\cong \Cok\delta'\oplus W''$, we have that $W''$ is a projective $(A',C')$-bimodule. So the composition $N''\to W'\twoheadrightarrow W''$ is a split surjective homomorphism. Hence we can write $N''\cong N'\oplus W''$ for some $(A',C')$-bimodule $N'$ and the sequence $(**)$ is isomorphic to a direct sum of $0 \to 0 \to W'' \xrightarrow{{\rm Id}} W'' \to 0$ and a sequence of the form
$$0\to Q'/\Ker\delta'\xrightarrow{\psi'}N'\to \Cok\delta'\to 0.$$
The map $f'$ in Theorem \ref{Ok1.2} (4) is the composition $Q'\xrightarrow{\phi'} Q'/\Ker\delta'\xrightarrow{\psi'} N'$.

\subsection{A descent criterion}\label{descentcriterion}

Keep the notation and assumptions in \S\ref{tilting}.  Assume that the field $k'$ above is finite. Let $k$ be a subfield of $k'$, and let $\Gamma:={\rm Gal}(k'/k)$. Assume that there are $k$-algebras $A$, $B$ such that $A'\cong k'\otimes_kA$, $B'\cong k'\otimes_k B$, respectively. 
For an $A'$-module
$U'$ and an automorphism $\sigma\in \Gamma$, denote by ${}^\sigma U'$ the $A'$-module which is equal to $U'$ as a module over the subalgebra $1\otimes A$ of $A'$, such that $x\otimes a$ acts on $U'$ as $\sigma^{-1}(x) \otimes a$ for all $a \in A$ and $x\in k'$.
 The $A'$-module $U'$ is {\it $\Gamma$-stable} if ${}^\sigma U'\cong U'$ for all $\sigma\in \Gamma$. $U'$ is said to be {\it defined over $k$}, if there is an $A$-module $U$ such that $U'\cong k'\otimes_k U$.
In this special case, $U'$ is $\Gamma$-stable, because for any $\sigma\in\Gamma$, the map sending $x\otimes u$ to $\sigma^{-1}(x)\otimes u$ is an isomorphism $k'\otimes_k U\cong {}^\sigma(k'\otimes_k U)$, where $u\in U$ and $x\in k'$.

In the proof of Theorem \ref{descenttool}, we will use \cite[Lemma 6.2]{Kessar_Linckelmann}.
The assumption in \cite[Lemma 6.2]{Kessar_Linckelmann} that $A/J(A)$ is separable is not needed (Kessar and Linckelmann approved this), because any finite-dimensional semisimple algebra over a finite field is separable (cf. \cite[page 130]{NT}).

\begin{theorem}\label{descenttool}
Keep the notation and assumptions above. Suppose that the set $\{T_i~|~i\in I_0\}$ is $\Gamma$-stable, i.e., for any $\sigma\in \Gamma$ and any $i\in I_0$, there exists $j\in I_0$ such that ${}^\sigma T_i\cong T_j$. Assume that the $(A',B')$-bimodule $M'$ is defined over $k$. The following holds.

\medskip\noindent{\rm(1)}. There exists a complex of $(A,B)$-bimodules $M(I_0)^{\bullet}$ such that $M'(I_0)^{\bullet}\cong k'\otimes_k M(I_0)^{\bullet}$ as complexes of $(A',B')$-bimodules.

\noindent{\rm(2)}. Let $C:={\rm End}_{K^b(A)}(M(I_0)^\bullet)^{\rm op}$, then $C'\cong k'\otimes_k C$ as $k'$-algebras.

\noindent {\rm (3)}. The right action of $B$ on $M(I_0)^\bullet$ induces a $k$-algebra monomorphism $\rho$ from $B$ to $C$, such that the diagram
$$\xymatrix{
  k'\otimes_k B \ar[r]^{{\rm Id}_{k'}\otimes \rho} \ar[d]_{\cong}
                & k'\otimes_k C \ar[d]^{\cong}  \\
        B' \ar[r]^{\rho'}       & C'             }$$
commutes. Moreover, $\rho$ is an isomorphism if and only if the monomorphism $\rho'$ in Theorem \ref{Ok1.2} (2) is an isomorphism.

\noindent{\rm(4)}. There exists a complex $N(I_0)^{\bullet}:=\cdots\to0\to Q\xrightarrow{f} N\to 0\to\cdots$, where $Q$ is a projective $(A,C)$-bimodule
of $(A,C)$-bimodules, such that $N'(I_0)^{\bullet}\cong k'\otimes_kN(I_0)^{\bullet}$ as complexes of $(A',C')$-bimodules, and $N(I_0)^{\bullet}$ induces a Rickard equivalence between $A$ and $C$.
\end{theorem}

\noindent{\it Proof.} (1). The method is inspired by the method, due to Kessar and Linckelmann, used in the proofs of \cite[Theorem 1.10]{Kessar_Linckelmann} and \cite[Theorem 1]{Huang}. By the assumptions, there is an $(A,B)$-bimodule $M$ such that $M'\cong k'\otimes_k M$. Since $M'$ induces a stable equivalence of Morita type between $A'$ and $B'$, it is projective as left $A'$-module and as right $B'$-module. By \cite[Lemma 4.4 (a)]{Kessar_Linckelmann}, $M$ is projective as left $A$-module and as right $B$-module. Since the set $\{T_i~|~i\in I_0\}$ is $\Gamma$-stable, it is easy to see that the set $\{P_i'\otimes_{k'} {Q'_i}^*|i\in I_0\}$ is $\Gamma$-stable. Indeed, for any $\sigma\in \Gamma$ and $i\in I_0$, there is $j\in I_0$ such that ${}^\sigma T_i\cong T_j$. Then ${}^\sigma (M'\otimes_{B'} T_i)\cong {}^\sigma M'\otimes_{B'} {}^\sigma T_j\cong M'\otimes_{B'} {}^\sigma T_j$. Since twisting by $\sigma$ is compatible with taking projective cover and taking dual, we have ${}^\sigma(P'_i\otimes_{k'}{Q'_i}^*)\cong P'_j\otimes_{k'}{Q'_j}^*$.
Hence the projective $(A',B')$-bimodule $P'(I_0)$ is $\Gamma$-stable. Then by \cite[Lemma 6.2 (c)]{Kessar_Linckelmann}, there is a projective $(A,B)$-bimodule $P(I_0)$ such that $P'(I_0)\cong k'\otimes_k P(I_0)$.

Next, we need to show that the homomorphism $\delta':=\delta'(I_0)$ can be chosen to be of the form ${\rm Id}_{k'}\otimes\delta$ for some bimodule homomorphism $\delta:=\delta(I_0):P(I_0)\to M$. Consider a projective cover $\pi: Z\to M$. Then $k'\otimes_k Z$ yields a projective cover of $M'$, hence $k'\otimes_k Z$ is isomorphic to the projective cover of $M'$ discussed in \S\ref{tilting}. By \cite[Lemma 6.2 (c)]{Kessar_Linckelmann} and the Noether-Deuring Theorem (see \cite[page 139]{CR}), $Z$ has a direct summand isomorphic to $P(I_0)$. So we just need to restrict $\pi$ to $P(I_0)$, and denote the map by $\delta(I_0)$, then $\delta:=\delta(I_0)$ is a desired map and
$$M(I_0)^\bullet:=\cdots\to 0\to P(I_0)\xrightarrow{\delta}M\to 0\to \cdots$$
is a desired complex.

\noindent(2). Note that as a complex of left $A$-modules (resp. $A'$-modules), every term of $M(I_0)^\bullet$ (resp. $M'(I_0)^\bullet$) is projective. Then by \cite[Proposition 3.5.43]{Zimmermann}, we have
$$C:={\rm End}_{K^b(A)}(M(I_0)^\bullet)^{\rm op}\cong {\rm End}_{D^b(A)}(M(I_0)^\bullet)^{\rm op}$$~~
and
$$C':={\rm End}_{K^b(A')}(M'(I_0)^\bullet)^{\rm op}\cong {\rm End}_{D^b(A')}(M'(I_0)^\bullet)^{\rm op}.$$
Now the statement follows from the Change of Ring Theorem (see \cite[Lemma 3.8.6]{Zimmermann}).

\noindent(3). The commutative diagram follows from the definition of the map $k'\otimes_k C\to C'$ and the definition of ``a module is defined over a subfield". By Theorem \ref{Ok1.2} (2), $\rho'$ is a monomorphism. Using the commutative diagram, we see that $\rho$ must be a monomorphism; if $\rho$ is surjective, then $\rho'$ is surjective; since $k'$ is a flat $k$-module, if $\rho$ is not surjective, then $\rho'$ could not be surjective.

\noindent(4). Since the algebras $A'$, $B'$, $C'$, and the homomorphism $\delta'$ are all defined over $k$, and since $k'$ is a flat $k$-module, the modules $\Ker\delta'$, $\I\delta'$ and $\Cok\delta'$ are defined over $k$. Now, it is a routine exercise to check that all modules and homomorphisms appeared in the construction of $f'$ (see \S\ref{f'}) are defined over $k$. In other words, all procedures in \S\ref{f'} can be realised over the field $k$. This implies the existence of $N(I_0)^{\bullet}$. By \cite[Lemma 4.4 (a)]{Kessar_Linckelmann}, $Q$ is a projective $(A,C)$-bimodule. By \cite[Proposition 4.5 (a)]{Kessar_Linckelmann}, $N(I_0)^{\bullet}$ induces a Rickard equivalence between $A$ and $C$.  $\hfill\square$

\section{Representation theory of $\SL$}\label{S2}

Let $G:=\SL$, and let $P:=\left\{\left[ {\begin{array}{*{20}{c}}
  1&b \\
  0&1
\end{array}} \right] \middle|~b\in \F_{p^n}\right\}$. Then $P$ is a Sylow $p$-subgroup of $G$, and $N_G(P)=\left\{\left[ {\begin{array}{*{20}{c}}
  a&b \\
  0&a^{-1}
\end{array}} \right] \middle|~a\in \F_{p^n}^\times,b\in \F_{p^n}\right\}$. We denote $N_G(P)$ by $H$.
We can see from \cite[\S2.1]{Oku} that simple $\F_{p^n}G$-modules and $\F_{p^n}H$-modules are absolutely simple, so $\F_{p^n}$ is a splitting field for both $G$ and $H$. Since every finite group has a finite splitting field, we assume that $k'$ is finite and take $k=\F_p$ in this section.  Since $\F_{p^n}$ is the minimal splitting field for $\SL$, we have $\F_{p^n}\subseteq k'$. We briefly review the simple modules of $k'G$ and $k'H$, after that we review some notation in \cite[\S2]{Oku} and \cite[\S 3]{Yo}.

Let $\sigma$ be the automorphism of the field $k'$, sending $x$ to $x^p$ for any $x\in k'$.
For $i\in \{1,\cdots, n\}$ and any $k'G$-module $U$, denote by $U^{(i)}$ the $k'G$-module ${}^{\sigma^i}U$.

\subsection{Simple modules of $k'G$ and $k'H$}\label{simplemodulesofSL}
For $\lambda\in\{0,1,\cdots, p-1\}$, Let $S_\lambda$ be the subspace of $k'[X,Y]$ consisting of all the homogeneous polynomials in indeterminates $X$ and $Y$ over $k'$ of degree $\lambda$. The group $G$ acts via invertible $k'$-linear transformations on $S_\lambda$: for any $g:=\left[ {\begin{array}{*{20}{c}}
  a&b \\
  c&d
\end{array}} \right]\in G$, and $f:=f(X,Y)\in S_\lambda$, $gf(X,Y):=f(aX+cY,bX+dY)$. For $\lambda\in\Lambda:=\{0,1,\cdots, p^n-1\}$ and its $p$-adic expansion $\lambda=\sum_{i=0}^{n-1}\lambda_ip^i$, it is well-known that
$$S_\lambda:=S_{\lambda_0}^{(0)}\otimes_{k'}S_{\lambda_1}^{(1)}\otimes_{k'}\cdots\otimes_{k'}S_{\lambda_{n-1}}^{(n-1)}$$
is a simple $k'G$-module, and $\{S_\lambda~|~\lambda\in\Lambda\}$ is a set of representatives of isomorphism classes of simple $k'G$-modules.

Let $\Lambda_0:=\Lambda-\{p^n-1\}$, $\mathcal{S}_1:=\{{\rm even~numbers~in}~\Lambda_0\}$, $\mathcal{S}_2:=\{{\rm odd~numbers~in}~\Lambda_0\}$. Denote the defect zero block, principal block, and the full defect non-principal block (if it exists) of $k'G$ by $b_0$, $b_1$, $b_2$ respectively.
Then $\{S_{p^n-1}\}$, $\{ S_\lambda~|~\lambda\in \mathcal{S}_1\}$, $\{ S_\lambda~|~\lambda\in \mathcal{S}_2\}$ are sets of representatives of isomorphism classes of simple $k'Gb_0$-, $k'Gb_1$-, $k'Gb_2$- modules, respectively if $p$ is odd and $\{S_{p^n-1}\}$, $\{S_\lambda~|~\lambda \in \Lambda_0\}$ are those of simple $k'Gb_0$-, $k'Gb_1$- modules, respectively if
$p=2$.

For $\lambda\in \Lambda_0$, let $T_\lambda$ be a $1$-dimensional vector space over $k'$ on which $h:=\left[ {\begin{array}{*{20}{c}}
  a&b \\
  0&a^{-1}
\end{array}} \right]\in H$ acts as scalar multiplication by $a^{\lambda}$. Then $\{T_\lambda~|~\lambda\in \Lambda_0\}$ is a set of representatives of isomorphism classes of simple $k'H$-modules. By definition, it is easy to check that ${\rm Soc}({\rm Res}_HS_\lambda)\cong T_\lambda$ for any $\lambda\in\Lambda_0$. Denote by $c_1$ the principal block of $k'H$ and $c_2$ the Brauer correspondent of $b_2$ (if $b_2$ exists). Then $\{ T_\lambda~|~\lambda\in \mathcal{S}_1\}$, $\{ T_\lambda~|~\lambda\in \mathcal{S}_2\}$ are sets of representatives of isomorphism classes of simple $k'Hc_1$-, $k'Hc_2$- modules, respectively if $p$ is odd and $\{T_\lambda~|~\lambda\in \Lambda_0\}$ is a set of representatives of isomorphism classes of simple $k'Hc_1$-modules if $p=2$.

\subsection{Some notation}

For $\lambda\in \Lambda_0$, define
$$\widetilde{\lambda}:=\left\{ \begin{gathered}
  0,~{\rm if}~\lambda=0; \hfill \\
  p^n-1-\lambda,~{\rm if}~\lambda\neq 0. \hfill \\
\end{gathered}  \right.$$
For a subset $\Omega\subseteq \Lambda_0$, let $\wt{\Omega}:=\{\wt{\lambda}~|~\lambda\in \Omega\}$. Clearly the map $\Lambda_0\to \Lambda_0, \lambda\mapsto \widetilde{\lambda}$ is a permutation on $\Lambda_0$ of order $2$.
For $\lambda$ and $\mu$ in $\Lambda_0$, define $\lambda\sim\mu$ if $S_\lambda\cong S_\mu^{(j)}$ for some integer $j$. Obviously ``$\sim$" is an equivalence relation on $\Lambda_0$. By using $p$-adic extensions, one easily verifies that $\lambda\sim \mu$ if and only if $\wt{\lambda}\sim\wt{\mu}$.

Let $I$ be $\mathcal{S}_1$ or $\mathcal{S}_2$ if $p$ is odd and be $\Lambda_0$ if $p=2$. Define the ordered equivalence classes with respect to ``$\sim$" as follows: let $J_{-1}$ and $\wt{J}_{-1}$ be empty sets (by convention), 
and $J_t$ the class containing the smallest $\lambda\in I-\cup_{u=-1}^{t-1}(J_u\cup \wt{J_u})$
for $t\geq 0$. Repeat this procedure until $t=s$ where $s$ satisfies $I=\cup_{u=-1}^{s}(J_u\cup \wt{J_u})$. Let $I_t:=\wt{J_t}$ and $K_t:=I_t\cup J_t$, so we have $I=\cup_{u=-1}^{s}K_u$.

\begin{remark}\label{gammastable}
{\rm Let $\Gamma:={\rm Gal}(k'/k)$ and let $\Gamma_0:={\rm Gal}(\F_{p^n}/k)$. We claim that for each $t\in\{0,\cdots,s\}$, the set $\{S_\lambda~|~\lambda\in I_t\}$ is $\Gamma$-stable. Since $\Gamma$ is a cyclic group, any subgroup of $\Gamma$ is a normal subgroup. By the fundamental theorem of Galois theory, any element $\rho\in\Gamma$ restricts to an automorphism $\tilde{\rho}$ of $\F_{p^n}$, and the map $\rho\mapsto \tilde{\rho}$ defines a surjective group homomorphism $\Gamma\to \Gamma_0$ with kernel ${\rm Gal}(k'/\F_{p^n})$. Since $\Gamma_0$ is a cyclic group of order $n$, generated by the automorphism $\F_{p^n}\to \F_{p^n}$, $x\mapsto x^p$, the set $\{\sigma^i~|~i=0,\cdots, n-1\}$ is a complete set of representatives for the pre-images of $\Gamma_0$ in $\Gamma$. Hence for every $\rho\in \Gamma$, $\rho=\sigma^i\rho_0$ for some $i$ and some $\rho_0\in {\rm Gal}(k'/\F_{p^n})$. For any $S_\lambda$ with $\lambda\in I_t$, we have
$${}^\rho S_\lambda={}^{\sigma^i}({}^{\rho_0}S_\lambda)\cong {}^{\sigma^i}S_\lambda=S_\lambda^{(i)},$$
where the second isomorphism holds because the $k'G$-module $S_\lambda$ is defined over $\F_{p^n}$.
By the choice of each $I_t$, $S_\lambda^{(i)}\in\{S_\lambda~|\lambda\in I_t\}$, as claimed.
Since ${\rm Soc}({\rm Res}_H S_\lambda)\cong T_\lambda$, the set $\{T_\lambda~|~\lambda\in I_t\}$ is also $\Gamma$-stable.
}
\end{remark}

\section{Rickard equivalences in $\SL$}\label{proof}

Keep the notation in Section \ref{S2}. In this section, we review the Rickard equivalences between blocks of $k'G$ and their Brauer correspondents constructed by Okuyama \cite{Oku} and Yoshii \cite{Yo}, and prove Theorem \ref{main2} for $G=\SL$. For defect zero blocks, Theorem \ref{main2} is trivial, so we only consider full defect blocks. Let $b$ be a full defect block of $k'G$, and let $c$ be the Brauer correspondent of $b$ in $k'H$. Let $A':=k'Gb$, $B':=k'Hc$, $A:=kGb$ and $B:=kHc$. Since $k = \F_p$, we note that $b \in kG$ and $c\in kH$ by Propositions \ref{blockidempotentofSL} and \ref{2.2}.

\begin{proposition}\label{lemma}
Multiplication by $b$ induces a unitary $k'$-algebra (resp. $k$-algebra) homomorphism $\rho'_0:B'\to A'$ (resp. $\rho_0:B\to A$), such that the diagram
$$\xymatrix{
   B \ar[r]^{\rho_0} \ar@{^{(}->}[d]
                & A \ar@{^{(}->}[d]  \\
        B' \ar[r]^{\rho'_0}       & ~A'             }$$
commutes. By the homomorphism $\rho'_0$ (resp. $\rho_0$), we can regard $A'$ (resp. $A$) as a left or right $B'$-module (resp. $B$-module). Hence the commutative diagram implies $A'\cong k'\otimes_k A$ as $(B',B')$-bimodules.
\end{proposition}

\noindent{\it Proof.} Consider $k'H$ as a $k'$-subalgebra of $k'G$. By the proof of \cite[Proposition 4.1.1 (b)]{Yo}, if the full defect non-principal block $b_2$ exists, then $b_2=c_2$ (as elements of $k'G$). Since $1=b_0+b_1+b_2=c_1+c_2$, we have $c_1=b_0+b_1$. So we have that $bc=b$ for $b\in \{b_1,b_2\}$ and its Brauer correspondent $c$. It follows that ``multiplication by $b$" defines a unitary $k'$-algebra homomorphism $\rho'_0:B'\to A'$, sending each $u\in B'$ to $bu$. Since all these block idempotents are contained in $A$ or $B$, the similar argument works for $A$ and $B$. By the construction of $\rho'_0$ and $\rho_0$, the diagram is obvious commutative. $\hfill\square$

\medskip The next proposition is summarized from \cite[Section 3]{Oku} and \cite[Proposition 4.1.1]{Yo}.

\begin{proposition}[see {\cite[Section 3]{Oku}, \cite[Proposition 4.1.1]{Yo}}]\label{OY}
 There is a sequence of $k'$-algebras $A'_0:=A',A'_1,\cdots,A'_s,A'_{s+1}$ such that the following holds for each $t\in \{0,1,\cdots,s\}$.

\medskip\noindent{\rm (1)}. $A'_t$ has a right $B'$-module structure. 

\noindent{\rm (2)}. The isomorphism classes of simple $A'_t$-modules are indexed by $I$.

\noindent{\rm (3)}. The algebras $A'_t$, $B'$ with the $(A'_t,B')$-bimodule $A'_t$ and the complex $A'_t(I_t)^\bullet$ of $(A'_t,B')$-bimodules satisfy the hypothesis in Theorem \ref{Ok1.2} ($A'_t$, $B'$, $A'_t$, $I$, $I_t$ instead of $A'$, $B'$, $M$, $I$, $I_0$, respectively). And $A'_{t+1}:=\End_{K^b(A'_t)}(A'_t(I_t)^\bullet)^{\rm op}$. Hence $A'_0,A'_1,\cdots,A'_{s+1}$ are all derived equivalent to $A'$.

\noindent{\rm (4)}. The right action of $B'$ on $A'_t(I_t)^\bullet$ induces a $k'$-algebra monomorphism $\rho'_{t+1}$ from $B'$ to $A'_{t+1}$. And $\rho'_{s+1}:B'\to A'_{s+1}$ is an isomorphism. For each $t\in\{0,\cdots,s\}$, the right $B'$-module structure of $A'_t$ in (1) is induced by $\rho'_t$ (here $\rho'_0$ is the homomorphism defined as in Proposition \ref{lemma}).

\noindent{\rm (5)}. For $t\in \{0,\cdots,s+1\}$, the $(A'_t,A')$-bimodule $A'_t\otimes_{B'} A'$ is isomorphic to a direct sum of a non-projective indecomposable module (denoted by $L'_t$) and a projective module.

\noindent{\rm (6)}.  Let $S_\lambda^t:=A'_t\otimes_{B'} T_\lambda$, for $\lambda\in \cup_{u=-1}^{t-1}K_u$, and let $S_\lambda^t:=L'_t\otimes_{A'} S_\lambda$, for $\lambda\in \cup_{u=t}^sK_u$. Then $\{S_\lambda^t~|~\lambda\in I\}$ is a set of representatives of isomorphism classes of simple $A'_t$-modules.

\noindent{\rm (7)} Denote the resulting Rickard complex for $A'_t$ and $A'_{t+1}$ in Theorem \ref{Ok1.2} (4) by
$$N'(I_t)^\bullet:=\cdots\to0\to Q'_t\xrightarrow{f'_t} N'_t\to 0\to\cdots.$$
Denote the $k'$-dual of $N'_t$ by ${N'_t}^*$. Then as $(A'_{t+1},B')$-bimodules, ${N'_t}^*\cong A'_{t+1}$; as $(A'_{t+1},A')$-bimodules, ${N'_t}^*\otimes_{A'_t}\otimes\cdots\otimes_{A'_1}{N'_0}^*$ is isomorphic to a direct sum of $L'_{t+1}$ and a projective $(A'_{t+1},A')$-bimodule. Here, $L'_{t+1}$ is defined as in (5). In particular, the complex
$$X'^\bullet:=N'(I_0)^\bullet\otimes_{A'_1}\cdots\otimes_{A'_s}N'(I_s)^\bullet$$
of $(A',B')$-bimodules  induces a Rickard equivalence between $A'$ and $B'$.
\end{proposition}

The next proposition is a slight refinement of Proposition \ref{OY} (7), which is suggested by an anonymous referee.

\begin{proposition}\label{referee}
Keep the notation of Proposition \ref{OY}. For each $t\in\{0,\cdots,s\}$, we have ${N'_t}^*\otimes_{A'_t}\otimes\cdots\otimes_{A'_1}{N'_0}^*\cong A'_{t+1}$ as $(A'_{t+1},B')$-bimodules.
The $(A'_{t+1},A')$-bimodule ${N'_t}^*\otimes_{A'_t}\otimes\cdots\otimes_{A'_1}{N'_0}^*$ is indecomposable and isomorphic to $L'_{t+1}$.
\end{proposition}

\noindent{\it Proof.} Since ${N'_t}^*\cong A'_{t+1}$ as $(A'_{t+1},B')$-bimodules for every $t\in\{0,\cdots,s\}$ (see Proposition \ref{OY} (7)), the first statement holds. Since $A'_0,A'_1,\cdots,A'_{s+1}$ are all derived equivalent to $A'$ and since $A'$ is an indecomposable $k'$-algebra, $A'_0,A'_1,\cdots,A'_{s+1}$ are indecomposable $k'$-algebras (see \cite[Lemma 6.7.12]{Zimmermann}). By Proposition \ref{OY} (3), for each $t\in \{0,\cdots,s\}$, the $(A'_t,B')$-bimodule $A'_t$ induces a stable equivalence between $A'_t$ and $B'$ and has no projective summands. Then by \cite[Theorem 2.1 (i)]{Lin96}, the $(A'_t,B)$-bimodule $A'_t$ is indecomposable. For $t=s+1$, we still have that the $(A'_t,B)$-bimodule $A'_t$ is indecomposable. So for each $t\in\{0,\cdots,s\}$, ${N'_t}^*\otimes_{A'_t}\otimes\cdots\otimes_{A'_1}{N'_0}^*\cong A'_{t+1}$ is an indecomposable non-projective $(A'_{t+1},B')$-bimodule, and hence it is also an indecomposable $(A'_{t+1},A')$-bimodule. Now by Proposition \ref{OY} (7), the second statement holds.  $\hfill\square$


\medskip The following is a descent of Proposition \ref{OY}.

\begin{proposition}\label{descent}
There are $k$-algebras $A_0:=A,A_1,\cdots,A_s,A_{s+1}$ such that for each $t\in \{0,1,\cdots,s\}$, the following holds.

\noindent{\rm (1)}. As $k'$-algebras, $A'_t\cong k'\otimes_kA_t$ and $A'_{s+1}\cong k'\otimes_k A_{s+1}$.

\noindent{\rm (2)}. There is an $(A_t,A)$-bimodule $L_t$, such that $L'_t\cong k'\otimes_k L_t$ as $(A'_t,A')$-bimodules.

\noindent{\rm (3)}. There is a complex $A_t(I_t)^\bullet$ of $(A_t,B)$-bimodules satisfying $A'_t(I_t)^\bullet\cong k'\otimes_k A_t(I_t)^\bullet$. The right action of $B$ on $A_t(I_t)^\bullet$ induces a $k$-algebra monomorphism $\rho_{t+1}$ from $B$ to $A_{t+1}$, such that the diagram
$$\xymatrix{
  k'\otimes_k B~ \ar[r]^{{\rm Id}_{k'}\otimes \rho_{t+1}~} \ar[d]_{\cong}
                & ~~k'\otimes_k A_{t+1} \ar[d]^{\cong}  \\
        B' \ar[r]^{\rho'_{t+1}}       & A'_{t+1}             }$$
commutes. Moreover, $\rho_{s+1}:B\to A_{s+1}$ is an isomorphism.

\noindent{\rm (4)}. There is a Rickard complexes
$$N(I_t)^\bullet:=\cdots\to0\to Q_t\xrightarrow{f_t} N_t\to 0\to\cdots$$
of $(A_t,A_{t+1})$-bimodules, such that $N'(I_t)^\bullet\cong k'\otimes_k N(I_t)^\bullet$ as complexes of $(A'_t,A'_{t+1})$-bimodules. In particular, the complex
$$X^\bullet:=N(I_0)^\bullet\otimes_{A_1}\cdots\otimes_{A_s}N(I_s)^\bullet$$
of $(A,B)$-bimodules  induces a Rickard equivalence between $A$ and $B$.
\end{proposition}

\noindent{\it Proof.} We proceed by induction on $t$. Since $A'_0=A'$ and $L'_0\cong A'$, let $A_0:=A$ and $L_0:=A$, then we have a $k'$-algebra isomorphism $A'_0\cong k'\otimes_k A_0$, and an $(A'_0,A')$-bimodule isomorphism $L'_0\cong k'\otimes_k L_0$. Hence (1),(2) hold for $t=0$.

Let $i\in \{0,1,\cdots, s\}$, assume that the statements (1), (2) hold for $t\leq i$, and (3), (4) hold for $t\leq i-1$ (when $i=0$, (3),(4) are empty propositions for $t\leq i-1$, hence hold). We are going to prove that (1), (2) hold for $t=i+1$, and (3), (4) hold for $t=i$. By the inductive hypothesis, $L'_i$ is defined over $k$, hence it is $\Gamma$-stable. Then by Remark \ref{gammastable}, the set $\{S_\lambda^i:=L'_i\otimes_{A'}S_\lambda~|~\lambda\in I_i\}$ is $\Gamma$-stable. Using the commutative diagram in (3) for $t=i-1$, we see that $A'_i$ is isomorphic to $k'\otimes_k A_i$ when they are regarded as $(A'_i,B')$-bimodules.
(When $i=0$, (3) is an empty proposition for $t=i-1$, but we still have an $(A'_0,B')$-bimodule isomorphism $A'_0\cong k'\otimes_k A_0$ by Proposition \ref{lemma}.)
So we can apply Theorem \ref{descenttool} to $A'_i$, $B'$, $A'_i$, $I$, $I_i$ instead of $A'$, $B'$, $M$, $I$, $I_0$, respectively. By Theorem \ref{descenttool} (1),(2), there is a complex $A_i(I_i)^\bullet$ of $(A_i,B)$-bimodules; let $A_{i+1}:=\End_{K^b(A_i)}(A_i(I_i)^\bullet)^{\rm op}$, then $A_{i+1}$ is a $k$-algebra and we have $A'_{i+1}\cong k'\otimes_k A_{i+1}$. Hence (1) holds for $t=i+1$.

The existence of the Rickard complex $N(I_{i})^\bullet$ follows by Theorem \ref{descenttool} (4), and thus (4) holds for $t=i$. Combining with the inductive hypothesis, now we have that the modules $N'_0,\cdots, N'_i$ in Proposition \ref{OY} (7) are all defined over $k$. By Proposition \ref{referee}, $L'_{i+1}$ is isomorphic to ${N'_i}^*\otimes_{A'_i}\cdots\otimes_{A'_1}{N'_0}^*$. Letting $L_{i+1}:={N_i}^*\otimes_{A_i}\cdots\otimes_{A_1}{N_0}^*$, then we have $L'_{i+1}\cong k'\otimes_k L_{i+1}$. Hence (2) holds for $t=i+1$.

The existence of the monomorphism $\rho_{i+1}$ and the commutative diagram in (3) for $t=i$ follow by Theorem \ref{descenttool} (3). Now we proved that the statements (1),(2),(3),(4) hold for all $t\in\{0,1,\cdots,s\}$.
Note that $\rho_{i+1}$ is an isomorphism if and only if $\rho'_{i+1}$ is an isomorphism. Hence $\rho_{s+1}$ is an isomorphism.   $\hfill\square$

\begin{lemma}\label{vertex}
The indecomposable $(A',B')$-bimodule $A'$ has a subgroup $Q$ of $\Delta P$ as a vertex.
\end{lemma}

\noindent{\it Proof.} The indecomposability of the $(A',B')$-bimodule $A'$ is by the proof of Proposition \ref{referee}. Since $A'$ is isomorphic to a direct summand of
$${\rm Res}_{G\times H^{\rm op}}^{G\times G^{\rm op}}{\rm Ind}_{\Delta P}^{G\times G^{\rm op}}(k'),$$
by the Mackey formula, $A'$ is isomorphic to a direct summand of
${\rm Ind}_{Q_x}^{G\times H^{\rm op}}(k')$,
where $Q_x=(G\times H^{\rm op})\cap x\Delta Px^{-1}$ for some $x\in G\times H^{\rm op}\setminus G\times G^{\rm op}/\Delta P$. Since $P$ is a trivial intersection subgroup of $G$, if $x\notin G\times H^{\rm op}$, then $Q_x=1$. So $x\in G\times H^{\rm op}$. We may choose $x=1$, and then $Q_x=\Delta P$. This completes the proof. $\hfill\square$


\begin{proposition}[Chuang]\label{splendid}
The Rickard complex $X^\bullet$ in Proposition \ref{descent} is splendid.
\end{proposition}

\noindent{\it Proof.} For each $t\in \{0,1,\cdots,s\}$, since $N(I_t)^\bullet$ is a Rickard complex for $A_t$ and $A_{t+1}$, by the definition of Rickard complex, $Q_t$ and $N_t$ are projective as left $A_t$-modules and as right $A_{t+1}$-modules. Moreover, by Theorem \ref{descenttool} (4), $Q_t$ is a projective $(A_t, A_{t+1})$-bimodule. By definition of tensor products of complexes, the degree 0 of $X^\bullet$ is the $(A,B)$-bimodule $N_0\otimes_{A_1}\cdots\otimes_{A_s}N_s$; all other terms of $X^\bullet$ is a direct sum of $(A,B)$-bimodules of the form $U_0\otimes_{A_1}\cdots\otimes_{A_s}U_s$, where $U_t\in \{Q_t,N_t\}$, and at least one $U_t$ should be $Q_t$. Hence $U_0\otimes_{A_1}\cdots\otimes_{A_s}U_s$ is projective as $(A,B)$-bimodule. By Proposition \ref{referee}, the $(A',B')$-bimodule $N'_0\otimes_{A'_1}\cdots\otimes_{A'_s}N'_s$ is isomorphic to ${L'_{s+1}}^*$. By Proposition \ref{OY} (5), the $(B',A')$-bimodule $L'_{s+1}$ is isomorphic to the $(B',A')$-bimodule $A'$. Since $A'$ is a symmetric algebra, ${L'_{s+1}}^*$ is isomorphic to the $(A',B')$-bimodule $A'$. By Lemma \ref{vertex}, ${L'_{s+1}}^*$ has a subgroup $Q\subseteq\{(u,u^{-1})|u\in P\}$ of $G\times H^{\rm op}$ as a vertex and has trivial source. Since ${L'_{s+1}}^*\cong k'\otimes_k(N_0\otimes_{A_1}\cdots\otimes_{A_s}N_s)$, by \cite[Lemma 5.1 and 5.2]{Kessar_Linckelmann}, $N_0\otimes_{A_1}\cdots\otimes_{A_s}N_s$ also has $Q$ as a vertex and has trivial source. Hence $X^\bullet$ is a splendid Rickard complex. $\hfill\square$

\medskip\noindent{\it Proof of Theorem \ref{main2} for $\SL$.} By Proposition \ref{blockidempotentofSL}, we may assume that $k=\F_p$. By Proposition \ref{descent} (4), there is a Rickard equivalence between $\F_p\SL$ and it Brauer correspondent algebra. By Proposition \ref{splendid}, the Rickard equivalence is splendid. $\hfill\square$

\section{On crossed products}

The construction of Rickard equivalences in $\GL$ relies on Marcus' theory on Rickard equivalences for group graded algebras. So we briefly review crossed products and some facts in this section. We refer to \cite[Definition 1.3.7 and 1.3.8]{Linckelmann} for the definitions of group graded algebras and crossed products.
We review some paragraphs in \cite[\S 1.3]{Linckelmann}.

Let $F$ be a field and $G$ a finite group. An $F$-algebra $B$ is called $G$-{\it graded} if, as an $F$-module, $B$ is a direct sum $B=\oplus_{x\in G} B_x$ satisfying $B_xB_y\subseteq B_{xy}$ for all $x,y\in G$. Note that the subspace $B_1$ indexed by the unit element of $G$ is a subalgebra of $B$. A {\it crossed product} of $A$ and $G$ is a $G$-graded $F$-algebra $B=\oplus_{x\in G} B_x$ such that $A=B_1$ and such that $B_x$ contains an invertible element in $B$, for all $x\in G$. For $x\in G$, choose $u_x \in B_x\cap B^\times$. Note that $u_x^{-1}\in B_{x^{-1}}$. For any $b\in B_x$ we have $b=(bu_x^{-1})u_x\in Au_x$, and hence $B_x=Au_x$. Conjugation by $u_x$ induces an algebra automorphism $\iota(u_x)$ on $A$, sending $a\in A$ to $u_xau_x^{-1}$. For any other choice $u_x'\in B_x\cap B^\times$ we have $u_x'=vu_x$ for some $v\in A^\times$, and hence $\iota(u_x')$ and $\iota(u_x)$ differ by an inner automorphism of $A$. By the discussion above, we have the following lemma.

\begin{lemma}\label{actiononsimplemodules}
Keep the notation above. $G$ acts on the set of isomorphism class of simple $A$-modules via the set of automorphisms $\{\iota(u_x)~|~x\in G\}$, in an obvious way.
\end{lemma}

The following lemmas are well-known, for the convenience of the reader, we sketch the proofs.

\begin{lemma}[Clifford theory for crossed products]\label{Clifford}
Let $B$ be a crossed product of an $F$-algebra $A$ and $G$. For any simple $B$-module $S$, the restriction ${\rm Res}_A^B(S)$ of $S$ to $A$ is a semisimple $A$-modules.
\end{lemma}

\noindent{\it Proof.} Mimic the proof of \cite[Theorem 1.9.9]{Linckelmann}. Let $T$ be a simple $A$-submodule of ${\rm Res}_A^B(S)$. For every $x\in G$, fix a $u_x\in B_x\cap B^\times$. It is easy to check that $u_xT$ is again a simple $A$-module and that $u_xu_yT=u_{xy}T$. It follows that the sum of all simple $A$-submodules of $S$ of the form $u_xT$, with $x\in G$, is a $B$-submodule of $S$. Since $S$ is simple, this implies that $S$ is the sum of the $u_xT$.  $\hfill\square$

\begin{lemma}[Maschke theory for crossed products]\label{Maschke}
Let $B$ be a crossed product of an $F$-algebra $A$ and $G$. Suppose that $|G|$ is invertible in $F$. Let $M$ be a $B$-module whose restriction to $A$ is semisimple as a $A$-module. Then $M$ is semisimple as a $B$-module.
\end{lemma}

\noindent{\it Proof.} Mimic the proof of \cite[Theorem 1.11.9]{Linckelmann}. Let $U$ be a $B$-submodule of $M$. We need to show that $U$ has a complement in $M$ as a $B$-module. Since $M$ is semisimple as an $A$-module, $U$ has a complement $V$ in $M$ as an $A$-module. For every $x\in G$, fix a $u_x\in B_x\cap B^\times$. Let $\pi:M\to U$ be the projection of $M$ onto $U$ with kernel $V$. Since
$V$ is an $A$-submodule of $M$, the map $\pi$ is an $A$-homomorphism but not necessarily a $B$-homomorphism. Define a map $\tau:M\to M$ by $\tau(m)=\frac{1}{{|G|}}\sum_{x\in G}u_x\pi(u_x^{-1}m)$ for all $m\in M$. Since $\pi$ is an $A$-homomorphism, the map $\tau$ does not depend on the choice of $u_x$. One checks that $\tau$ is a projection of $M$ to $U$ as a $B$-module, and hence ${\rm ker}(\tau)$ is a complement of $U$ in $M$. Thus $M$ is semisimple. $\hfill\square$

\begin{lemma}\label{radical}
Let $B$ be a crossed product of an $F$-algebra $A$ and $G$. Suppose that $|G|$ is invertible in $F$. Then $J(B)=J(A)B=BJ(A)$.
\end{lemma}

\noindent{\it Proof.} Mimic the proof of \cite[Theorem 1.11.10]{Linckelmann}. By Lemma \ref{Clifford}, every simple $B$-module restricts to a semisimple $A$-module, hence is annihilated by $J(B)$. Thus $J(A)\subseteq J(B)$. For every $x\in G$, fix a $u_x\in B_x\cap B^\times$. Since each $u_x$ induces an automorphism of $A$, hence stabilises $J(A)$, and therefore $J(A)B=BJ(A)$ is an ideal contained in $J(B)$. In order to show that $J(B)\subseteq J(A)B$, it suffices to show that $B/J(A)B$ is semisimple as a $B$-module. Since $|G|$ is invertible in $F$, it suffices, by Lemma \ref{Maschke}, to show that $B/J(A)B$ is semisimple as an $A$-module. This is clear since $J(A)$ annihilates $B/J(A)B$. $\hfill\square$

\begin{proposition}\label{Moritaeq}
Let $B$ and $C$ be crossed products of an $F$-algebra $A$ and $G$. Then $B$ and $C$ are graded Morita equivalent in the sense of \cite[Definition 3.3]{Marcuscomm}.
\end{proposition}

\noindent{\it Proof.} Denote by $\Delta$ the diagonal subalgebra
$$\Delta(B,C):=\oplus_{x\in G} B_x\otimes_FC_{x^{-1}}$$
of $B\otimes_FC^{\rm op}$. For every $x\in G$, fix a $u_x\in B_x\cap B^\times$ and a $v_x\in C_x\cap C^\times$. Then we have that $B_x=u_xA$ and $C_x=Av_x$. For any $\delta_x:=u_xa_1\otimes a_2v_{x^{-1}}\in B_x\otimes_FC_{x^{-1}}$, where $a_1,a_2\in A$, and any $a\in A$, define $\delta_x\cdot a:=u_xa_1aa_2v_{x^{-1}}$.
It is obvious that with this action, the $(A,A)$-bimodule $A$ extends to a $\Delta$-module. Since the $(A,A)$-bimodule $A$ induces a Morita self-equivalence of the $F$-algebra $A$, by \cite[Theorem 3.4]{Marcuscomm}, the $(B,C)$-bimodule $B\otimes_A A$ and the $(C,B)$-bimodule $C\otimes_A A$ induce a graded Morita equivalence between $B$ and $C$.  $\hfill\square$

\section{Rickard equivalences in ${\rm GL}(2,p^n)$}\label{s6}

Let $G:=\SL$, $P:=\left\{\left[ {\begin{array}{*{20}{c}}
  1&b \\
  0&1
\end{array}} \right] \middle|~b\in \F_{p^n}\right\}$, $H:=N_G(P)$, and $\tilde{G}:=\GL$. By calculating orders, we see that $P$ is also a Sylow $p$-subgroup of $\tilde{G}$. Set $\tilde{H}:=N_{\tilde{G}}(P)$. By Frattini argument, $\tilde{G}=G\tilde{H}$, so we have $\tilde{G}/G\cong \tilde{H}/H\cong C_{p^n-1}$, where $C_{p^n-1}$ denotes a cyclic group of order $p^n-1$. Since every finite group has a finite splitting field, we assume that $k'$ is finite and take $k=\F_p$ in this section. Let $b$ be a full defect block of $k'G$. By Lemma \ref{blockidempotentofGL}, we know that $b$ is the sum of all blocks of $k'\tilde{G}$ which covers $b$. Let $c$ be the Brauer correspondent of $b$ in $k'H$. By the same argument, $c$ is the sum of all blocks of $k'\tilde{H}$ which covers $c$. Our task in this section is to prove the following proposition.

\begin{proposition}\label{blocksumofGL}
$k\tilde{G}b$ and $k\tilde{H}c$ are splendidly Rickard equivalent.
\end{proposition}

When $b$ is the principal block of $k'G$, Marcus \cite{Marcus03} proved that $k'\tilde{G}b$ and $k'\tilde{H}c$ are Rickard equivalent by showing that Okuyama's equivalences are compatible with $p'$-extensions. The main tool in Marcus' proof is \cite[Proposition 3.13]{Marcus03}. We first review \cite[Proposition 3.13]{Marcus03} and then give a descent criterion for \cite[Proposition 3.13]{Marcus03}.

Let $A'$ and $B'$ be split finite-dimensional symmetric $k'$-algebras having no semisimple summand. Let $M'$, $T_i$, $Q'_i$, $P'_i$, $\delta'_i$, $I$, $I_0$, $P'(I_0)$, $\delta':=\delta'(I_0)$, $M'(I_0)^\bullet$ and $N'(I_0)^\bullet$ be as in \S\ref{tilting}. Assume that $R'$ (resp. $S'$) is a crossed product graded by a finite group $T$ with $R'_1=A'$ (resp. $S'_1=B'$). Then $I$ can be regarded as a $T$-set via the action of $T$ on the set $\{T_i|i\in I\}$ of simple $B'$-modules (see Lemma \ref{actiononsimplemodules}). Denote by $\Delta'$ the diagonal subalgebra
$$\Delta(R',S'):=\oplus_{t\in T} R'_t\otimes_{k'}S'_{t^{-1}}$$
of $R'\otimes_{k'}S^{'{\rm op}}$. Clearly $\Delta'$ is a crossed product of $\Delta'_1:=A'\otimes_{k'}{B'}^{\rm op}$ and $T$.

Write $C':={\rm End}_{K^b(A')}(M'(I_0)^\bullet)^{\rm op}$ and $E':={\rm End}_{K^b(R')}(R'\otimes_{A'}M'(I_0)^\bullet)^{\rm op}$. The complex $R'\otimes_{A'}M'(I_0)^\bullet$ is $T$-graded, with $1$-component $M'(I_0)^\bullet$. Hence by \cite[Lemma 1.7 (a)]{Marcus03}, $E'$ is a $T$-graded algebra with $1$-component
$$E'_1\cong \End_{K^b(R'\mbox{-}{\rm Gr})}(R'\otimes_{A'}M'(I_0)^\bullet)^{\rm op}\cong {\rm End}_{K^b(A')}(M'(I_0)^\bullet)^{\rm op}=C'.$$
Here, $R'\mbox{-}{\rm Gr}$ denotes the category of $G$-graded $R'$-modules; the second isomorphism holds because a grade-preserving $R'$-homomorphism is determined by its restriction on the $1$-component. Since the $1$-component of $E'$ is isomorphic to $C'$, $C'$ can be regarded as a $\Delta(E',E')$-module.

\begin{proposition}[{\cite[Proposition 3.13]{Marcus03}}]\label{Marcusprop}
Assume that $T$ is a $p'$-group, $M'$ is a $\Delta'$-module and $I_0$ is a $T$-subset of $I$. Then the following holds.

\noindent{\rm (a)}. $M'(I_0)^\bullet$ extends to a complex of $\Delta'$-modules.

\noindent{\rm (b)}. $E'$ is a crossed product, and there is a graded stable equivalence of Morita type between $E'$ and $S'$.

\noindent{\rm (c)}. $N'(I_0)^\bullet$ extends to a complex of $\Delta(R',E')$-modules.
\end{proposition}

As noted in the last paragraph of \cite[page 192]{Marcus03}, all simple $k'G$-modules (resp. $k'H$-modules) are $C_{p^n-1}$-invariant since they extend to $\tilde{G}$ (resp. $\tilde{H}$).  So Okuyama and Yoshii's equivalences (Proposition \ref{OY} (7)), Proposition \ref{Marcusprop} (c) and \cite[Theorem 4.8]{Marcuscomm} imply that there is a Rickard equivalence between $k'\tilde{G}b$ and $k'\tilde{H}c$.

Assume that there are $k$-algebras $A$, $B$ such that $A'\cong k'\otimes_k A$, $B'\cong k'\otimes_k B$. Assume that $R$ (resp. $S$) is a crossed product graded by the finite group $T$ with $R_1=A$ (resp. $S_1=B$) satisfying $R'\cong k'\otimes_k R$ (resp. $S'\cong k'\otimes_k S$) as crossed products. That means, for each $t\in T$, we have $R'_t\cong k'\otimes_k R_t$ (resp. $S'_t\cong k'\otimes_k S_t$). Denote by $\Delta$ the diagonal subalgebra
$$\Delta(R,S):=\oplus_{t\in T} R_t\otimes_kS_{t^{-1}}$$
of $R\otimes_k S^{\rm op}$. Clearly we have $\Delta'\cong k'\otimes_k \Delta$. Let $\Gamma:={\rm Gal}(k'/k)$. Suppose that the set $\{T_i|i\in I_0\}$ is $\Gamma$-stable and that there is an $(A,B)$-bimodule $M$ satisfying $M'\cong k'\otimes_k M$. Then by Theorem \ref{descenttool}, there exists a complex of $(A,B)$-bimodules $M(I_0)^{\bullet}$ such that $M'(I_0)^{\bullet}\cong k'\otimes_k M(I_0)^{\bullet}$;
let $C:={\rm End}_{K^b(A)}(M(I_0)^\bullet)^{\rm op}$, then $C'\cong k'\otimes_k C$ as $k'$-algebras; there exists a complex $N(I_0)^{\bullet}$ of $(A,C)$-bimodules such that $N'(I_0)^{\bullet}\cong k'\otimes_kN(I_0)^{\bullet}$, and $N(I_0)^{\bullet}$ induces a Rickard equivalence between $A$ and $C$. Write  $E:={\rm End}_{K^b(R)}(R\otimes_{A}M(I_0)^\bullet)^{\rm op}$.

\medskip The following is a descent of Proposition \ref{Marcusprop}.

\begin{proposition}\label{Marcusdescent}
Keep the notation above. Assume that $T$ is a $p'$-group, $I_0$ is a $T$-subset of $I$, $M$ is a $\Delta$-module and $M'\cong k'\otimes_k M$ as $\Delta'$-modules. Then the following holds.

\noindent{\rm (a)}. $M(I_0)^\bullet$ extends to a complex of $\Delta$-modules.

\noindent{\rm (b)}. $E$ is a crossed product, and there is a graded stable equivalence of Morita type between $E$ and $S$.

\noindent{\rm (c)}. $N(I_0)^\bullet$ extends to a complex of $\Delta(R,E)$-modules. Moreover, $N'(I_0)^\bullet\cong k'\otimes_k N(I_0)^\bullet$ as complexes of $\Delta(R',E')$-modules.
\end{proposition}

\noindent{\it Proof.} (a). In \S\ref{tilting}, we see that the $(A',B')$-bimodule $\oplus_{i\in I}P'_i\otimes_{k'} {Q'}_i^*$ (let us denote it by $P'(M')$) is a projective cover of ${}_{A'}M'_{B'}$. Let $P'$ be a projective cover of $M'$ in the category of $\Delta'$-modules. Since $\Delta'$ is a crossed product of $\Delta'_1:=A'\otimes_{k'}{B'}^{\rm op}$ and $T$, and $T$ is a $p'$-group, by Lemma \ref{radical}, we have $P'/J(\Delta')P'=P'/J(\Delta_1')P'$. So ${\rm Res}_{\Delta_1'}^{\Delta'}(P')$ is a projective cover of ${}_{A'}M'_{B'}$, hence isomorphic to $P'(M')$. Thus we can assume that ${\rm Res}_{\Delta_1'}^{\Delta'}(P')=P'(M')$.

Let $P$ be a projective cover of $M$ in the category of $\Delta$-modules. Then $k'\otimes_k P$ yields a projective cover of $M'$, hence $k'\otimes_k P \cong P'$. Denote by $\pi'$ (resp. $\pi$) the canonical surjection $P'\twoheadrightarrow M'$ (resp. $P\twoheadrightarrow M$), we have the following commutative diagram:
$$\xymatrix{
  k'\otimes_k P~ \ar@{->>}[r]^{{\rm Id}_{k'}\otimes \pi} \ar[d]_{\cong}
                & ~~k'\otimes_k M \ar[d]^{\cong}  \\
        P' \ar@{->>}[r]^{\pi'}       & M'             }$$
By a similar argument as in the previous paragraph, the $(A,B)$-bimodule $P(M):={\rm Res}_{\Delta_1}^\Delta(P)$ is a projective cover of ${}_AM_B$.
So we have $k'\otimes_k P(M)\cong P'(M')$ as $(A',B')$-bimodules.

By the proof of Proposition \ref{Marcusprop} (a) (see \cite{Marcus03}), the $\Delta_1'$-summand $P'(I_0)$ of $P'(M)$ is a $T$-invariant $\Delta_1'$-module, hence $P'(I_0)$ extends to a $\Delta'$-summand of $P'$. Recall that $M'(I_0)^\bullet$ is a complex of the form $\cdots\to 0\to P'(I_0)\xrightarrow{\delta'}M'\to 0\to \cdots$. By the definition of $\delta'$, we see that as a map, $\delta'$ is the restriction of $\pi'$ to $P'(I_0)$. In other words, $\delta'$ extends to the $\Delta'$-homomorphism $\pi'$. In the proof of Theorem \ref{descenttool} (1), we showed that the $(A,B)$-bimodule $P(M)$ has a direct summand isomorphic to $P(I_0)$; let $\delta:=\pi|_{P(I_0)}$ (consider it as a homomorphism of $(A,B)$-bimodule), then the complex $M(I_0)^\bullet$ is of the form
$$\cdots\to 0\to P(I_0)\xrightarrow{\delta}M\to 0\to \cdots.$$ Since $P'(I_0)\cong k'\otimes_k P(I_0)$ as $(A',B')$-bimodules and since $P'(I_0)$ is $T$-invariant, we can deduce that $P(I_0)$ should also be $T$-invariant. Hence $P(I_0)$ extends to a $\Delta$-summand of $P$, and $\delta$ is also a $\Delta$-homomorphism (because $\delta=\pi|_{P(I_0)}$, and $\pi|_{P(I_0)}$ is a $\Delta$-homomorphism).

\medskip\noindent(b). By \cite[Lemma 2.6]{Marcuscomm}, $R\otimes_A M(I_0)^\bullet$ is a complex of $T$-graded $(R,S)$-bimodules, so the right multiplication gives a map $S\to E$ of $T$-graded algebras. Since $S$ is a crossed product, each component of $S$ contains an invertible element in $S$. Hence each component of $E$ contains an invertible element in $E$, which implies that $E$ is a crossed product.

By a similar argument as in the paragraph preceding Proposition \ref{Marcusprop}, we see that the $1$-component of $E$ is isomorphic to $C$, so $C$ can be regarded as a $\Delta(E,E)$-module. Via the $T$-graded algebra homomorphism $S\to E$, $C$ can also be regarded as a $\Delta(S,E)$-module, a $\Delta(E,S)$-module, or a $\Delta(S,S)$-module.

Since $R'\cong k'\otimes_k R$, $S'\cong k'\otimes_k S$, $A'\cong k'\otimes_kA$, $B'\cong k'\otimes_k B$, and $M'(I_0)^\bullet\cong k'\otimes_k M(I_0)^\bullet$, we see that $R'\otimes_{A'} M'(I_0)^\bullet\cong k'\otimes_k (R\otimes_AM(I_0)^\bullet)$ as complexes of $T$-graded $(R',S')$-bimodules. Then by the Change of Ring Theorem (see \cite[Lemma 3.8.6]{Zimmermann}), we have $E'\cong k'\otimes_k E$. By the proof of Proposition \ref{Marcusprop} (b) (see \cite{Marcus03}), we have that
\begin{equation}\label{7.1}
C'\cong B'\oplus P'_1~~{\rm as}~~\Delta(S',S')\mbox{-modules},\tag{7.3.1}
\end{equation}
\begin{equation}
C'\otimes_{B'}C'\cong C'\oplus P'_2~~{\rm as}~~\Delta(E',E')\mbox{-modules},~~~~~~~~~\tag{7.3.2}
\end{equation}
where $P'_1$ is a projective $\Delta(S',S')$-module, and $P'_2$ is a projective $\Delta(E',E')$-module.

Since the complex $R'\otimes_{A'} M'(I_0)^\bullet$ of $T$-graded $(R',S')$-modules is defined over $k$, we can deduce that
$C'\cong k'\otimes_k C$ as $\Delta(E',S')$-modules, as $\Delta(S',E')$-modules, as $(E',E')$-modules, and also as $\Delta(S',S')$-modules.
By the assumption on $B'$, we also have $B'\cong k'\otimes_k B$ as $\Delta(S',S')$-modules. It follows that
\begin{equation}\label{7.3}
{}^\sigma C'\cong C',~~~{}^\sigma B'\cong B \tag{7.3.3}
\end{equation}
as $\Delta(S',S')$-modules, for any $\sigma\in \Gamma$. By (\ref{7.1}), (\ref{7.3}) and Krull-Schmidt Theorem, we have ${}^\sigma P'_1\cong P'_1$ as $\Delta(S',S')$-modules. By \cite[Lemma 6.2 (c)]{Kessar_Linckelmann}, there exists a projective $\Delta(S,S)$-module $P_1$ satisfying $P'_1\cong k'\otimes_k P_1$. Then by the Noether-Deuring Theorem (see \cite[page 139]{CR}) we have that
$$C\cong B\oplus P_1~~{\rm as}~~\Delta(S,S)\mbox{-modules}.$$
A similar argument shows that there exists a $\Delta(E,E)$-module $P_2$, such that
$$C\otimes_BC\cong C\oplus P_2~~{\rm as}~~\Delta(E,E)\mbox{-modules}.~~~~~~~~$$
Using \cite[Theorem 5.4]{Marcuscomm}, we obtain that the $G$-graded $(E,S)$-bimodule $E$ induces a graded stable equivalence of Morita type between $E$ and $S$.

\medskip\noindent(c). Since we have proved (a), (b), and since we showed in \S\ref{descentcriterion} that all modules and complexes appeared in the proof of Proposition \ref{Marcusprop} (c) (see \cite{Marcus03}) are defined over $k$, the same argument as in the proof of Proposition \ref{Marcusprop} (c) shows that $N(I_0)^\bullet$ extends to a complex of $\Delta(R,E)$-modules. Then $k'\otimes_k N(I_0)^\bullet$ yields a $\Delta(R',E')$-module which is isomorphic to $N'(I_0)$, hence (c) holds.  $\hfill\square$

\medskip\noindent {\it Proof of Proposition \ref{blocksumofGL}.}  We adopt the notation in the first paragraph of this section. Let $R':=k'\tilde{G}b$, $S':=k'\tilde{H}c$. Then $R'$ (resp. $S'$) is a crossed product graded by the cyclic group $C_{p^n-1}\cong \tilde{G}/G$, with $R'_1=k'Gb$ (resp. $S'_1=k'Hc$). We borrow the notation in Proposition \ref{OY}.  Recall that for $t\in\{0,1,\cdots,s\}$, we have a $k'$-algebra $A'_{t+1}={\rm End}_{K^b(A'_t)}(A'_t(I_t)^\bullet)^{\rm op}$ (see Proposition \ref{OY} (3)). Set $E'_0=R'$, and define $E'_{t+1}:={\rm End}_{K^b(E'_t)}(E'_t\otimes_{A'_t}A_t'(I_t)^\bullet)^{\rm op}$ inductively for each $t\in \{0,1,\cdots,s\}$.
For $t\in\{0,1,\cdots,s\}$ denote by $\Delta'_t$ the diagonal subalgebra
$\Delta(E'_t,S')$.
By Proposition \ref{Marcusprop}, the complex $N'(I_t)^\bullet$
of $(A'_t,A'_{t+1})$-bimodules extends to a complex of $\Delta(E'_t,E'_{t+1})$-modules. Hence by \cite[Lemma 2.9]{Marcuscomm}, the complex $N'(I_0)^\bullet\otimes_{A'_1}\cdots\otimes_{A'_s}N'(I_s)^\bullet$ of $(A'_0,A'_{s+1})$-bimodules extends to a complex of $\Delta(E'_0,E'_{s+1})$-modules. Since both $E'_{s+1}$ and $S'$ are crossed products of $B'$ and $T$, by the proof of Proposition \ref{Moritaeq}, the $(B',B')$-bimodule $B'$ extends to a $\Delta(E'_{s+1},S')$-module. If we identify $A'_{s+1}$ and $B'$ via the isomorphism in Proposition \ref{OY} (4), then the complex $X'^\bullet:=N'(I_0)^\bullet\otimes_{A'_1}\cdots\otimes_{A'_s}N'(I_s)^\bullet$ of $(A',B')$-bimodules extends to a complex of $\Delta(R',S')$-modules. Using \cite[Corollary 3.9 (b)]{Marcusonequivalences}, we see that the complex $\widetilde{X}'^\bullet:=(R'\otimes_{k'} S'^{\rm op})\otimes_{\Delta(R',S')} X'^\bullet$ is a splendid Rickard complex for $R'$ and $S'$.

Actually, in the paragraph above, we reviewed Marcus' method for proving that $R'$ and $S'$ are Rickard equivalent, and we proved further that $R'$ and $S'$ are splendidly Rickard equivalent.

Let $R:=k\tilde{G}b$, $S:=k\tilde{H}c$. Then $R$ (resp. $S$) is a crossed product graded by the cyclic group $C_{p^n-1}\cong \tilde{G}/G$, with $R_1=kGb$ (resp. $S_1=kHc$). We borrow the notation in Proposition \ref{descent}. Recall that for $t\in\{0,1,\cdots,s\}$, we have a $k$-algebra $A_{t+1}={\rm End}_{K^b(A_t)}(A_t(I_t)^\bullet)^{\rm op}$ (see the proof of Proposition \ref{descent}). Set $E_0=R$, and define $E_{t+1}:={\rm End}_{K^b(E_t)}(E_t\otimes_{A_t}A_t(I_t)^\bullet)^{\rm op}$ inductively for each $t\in \{0,1,\cdots,s\}$.
For $t\in\{0,1,\cdots,s\}$, denote by $\Delta_t$ the diagonal subalgebra
$\Delta(E_t,S)$. Using the Change of Ring Theorem, it is easy to see that as $T$-graded algebras $E'_i\cong k'\otimes_k E_i$, so we have $\Delta'_t\cong k'\otimes_k \Delta_t$. To prove Proposition \ref{blocksumofGL}, we use Proposition \ref{Marcusdescent} to show that the complex $N(I_t)^\bullet$
of $(A_t,A_{t+1})$-bimodules extends to a complex of $\Delta(E_t,E_{t+1})$-modules. The assumption of Proposition \ref{Marcusdescent} requires us to show that the following two conditions hold:

\noindent(i). the $(A_t,B)$-bimodule $A_t$ extends to a $\Delta_t$-module;

\noindent(ii). $A'_{t}\cong k'\otimes_k A_t$ as $\Delta'_t$-modules.

\noindent Note that we had proved these two conditions in the proof of Proposition \ref{Marcusdescent} (b) ($A'_{t}$, $E'_t$, $\Delta'_t$, $A_t$, $E_t$, $\Delta_t$ instead of $C'$, $E'$, $\Delta(E',S')$, $C$, $E$, $\Delta(E,S)$, respectively). So by Proposition \ref{Marcusdescent}, the complex $N(I_t)^\bullet$
of $(A_t,A_{t+1})$-bimodules extends to a complex of $\Delta(E_t,E_{t+1})$-modules. Now, a similar argument as in the first paragraph shows that the complex $X^\bullet:=N(I_0)^\bullet\otimes_{A_1}\cdots\otimes_{A_s}N(I_s)^\bullet$ of $(A,B)$-bimodules extends to a complex of $\Delta(R,S)$-modules, and the complex $\widetilde{X}^\bullet:=(R\otimes_k S^{\rm op})\otimes_{\Delta(R,S)} X^\bullet$ is a splendid Rickard complex for $R$ and $S$.  $\hfill\square$

\section{Proof of Theorem \ref{main2} for $\GL$}\label{proofforGL}

Let $G:=\SL$, $P:=\left\{\left[ {\begin{array}{*{20}{c}}
  1&b \\
  0&1
\end{array}} \right] \middle|~b\in \F_{p^n}\right\}$, $H:=N_G(P)$, $\tilde{G}:=\GL$, and $\tilde{H}:=N_{\tilde{G}}(P)$. Let $b$ be a block of $k'G$. By Lemma \ref{blockidempotentofGL}, we know that $b$ is the sum of all blocks of $k'\tilde{G}$ which covers $b$; conversely, every block of $k'\tilde{G}$ is either full defect or defect zero, and covers a unique block of $k'G$. To prove Theorem \ref{main2} for $\GL$, we only need to consider full defect blocks. So we assume that $b$ is of full defect.
Let $c$ be the Brauer correspondent of $b$ in $k'H$. By the same argument, $c$ is the sum of all blocks of $k'\tilde{H}$ which covers $c$. Assume that $b=b_1+\cdots+b_r$, where $b_1,\cdots,b_r$ are blocks of $k'\tilde{G}$. Let $c_i$ ($i=1,\cdots,r$) be a block of $k'\tilde{H}$, corresponding to $b_i$ via the Brauer correspondence. It is easy to see that $c=c_1+\cdots+c_r$. Since $b_1$ can be any full defect block of $k'G$, so it suffices to prove Theorem \ref{main2} for the block $b_1$. From now on we assume that $k=\F_p[b_1]$. By Proposition \ref{2.2}, we have $c_1\in k\tilde{H}$. Now Theorem \ref{main2} for $\GL$ is equivalent to the following statement.

\begin{proposition}\label{8.1}
$k\tilde{G}b_1$ and $k\tilde{H}c_1$ are splendidly Rickard equivalent.
\end{proposition}

The following lemma is well-known, see e.g. \cite[Lemma 1.9]{Ruh} for a proof.

\begin{lemma}\label{summand}
Let $F$ be a field, and let $G_1$, $G_2$ be finite groups. Let $e$ (resp. $f$) be an idempotents in the center of $FG_1$ (resp. $FG_2$). Assume that a complex $C$ of $(FG_1e,FG_2f)$-bimodules induces a Rickard equivalence between $FG_1e$ and $FG_2f$. Let $e=e_1+\cdots+e_r$ be a decomposition of $e$ into blocks of $FG_1$. Then for each $i\in\{1,\cdots, r\}$, there exists a unique block $f_i$ of $FG_2$ such that $e_iCf_i$ is not homotopy equivalent to $0$. Moreover, $e_iCf_i$ induce a Rickard equivalence between $FG_1e_i$ and $FG_2f_i$.
\end{lemma}

By Proposition \ref{blocksumofGL}, $\F_p\tilde{G}b$ and $\F_p\tilde{H}c$ are splendidly Rickard equivalent. Since
$\F_p\subseteq k=\F_p[b_1]$, there is a complex $X$ of $(k\tilde{G}b,k\tilde{H}c)$-bimodules inducing a splendid Rickard between $k\tilde{G}b$ and $k\tilde{H}c$. Let $X':=k'\otimes_k X$, then $X'$ is a complex of $(k'\tilde{G}b,k'\tilde{H}c)$-bimodules inducing a splendid Rickard between $k'\tilde{G}b$ and $k'\tilde{H}c$.
By Lemma \ref{summand}, there is a unique $c_t$, where $t\in \{1,\cdots, r\}$, such that $b_1X'c_t$ is not homotopy equivalent to $0$, and $b_1X'c_t$ induces a splendid Rickard equivalence between $k'\tilde{G}b_1$ and $k'\tilde{H}c_t$. We don't know whether $c_t$ is $c_1$, the Brauer correspondent of $b_1$.

\begin{lemma}\label{filedvalue}
Keep the notation above, we have $k=\F_p[b_1]=\F_p[c_t]$, and $b_1Xc_t$ induces a splendid Rickard equivalence between $k\tilde{G}b_1$ and $k\tilde{H}c_t$.
\end{lemma}

\noindent{\it Proof.} Since every finite group has a finite splitting field, we may assume that $k'$ is finite.
For any $\sigma\in {\rm Gal}(k'/k)$, $\sigma(c_t)$ is a block of $k'\tilde{H}$. Since $b_1X'c_t$ is not homotopy equivalent to $0$,  $b_1X'\sigma(c_t)\cong b_1({}^\sigma X')\sigma(c_t)$ is not homotopy equivalent to $0$. The uniqueness of $c_t$ forces $c_t=\sigma(c_t)$, so $c_t\in k\tilde{H}$. By the same argument, we can deduce that the coefficients of $b_1$ are contained in the field $\F_p[c_t]$, so $k=\F_p[c_t]$. The second statement follows by \cite[Proposition 4.5 (a)]{Kessar_Linckelmann} (or using Lemma \ref{summand} again) $\hfill\square$

\medskip By \cite[Lemma 10.2.6]{Rou98}, if we can prove the following proposition, then Proposition \ref{8.1} holds.

\begin{proposition}
$k\tilde{H}c_t$ and $k\tilde{H}c_1$ are splendidly Morita equivalent.
\end{proposition}

\noindent{\it Proof.} We first prove that $k'\tilde{H}c_t$ and $k'\tilde{H}c_1$ are splendidly Morita equivalent.
Recall that $$\tilde{H}=\left\{\left[ {\begin{array}{*{20}{c}}
  a&b \\
  0&d
\end{array}} \right] \middle|~a,d\in \F_{p^n}^\times,~b\in \F_{p^n}\right\}~~~{\rm and}~~~H=\left\{\left[ {\begin{array}{*{20}{c}}
  a&b \\
  0&a^{-1}
\end{array}} \right] \middle|~a\in \F_{p^n}^\times,~b\in \F_{p^n}\right\}.$$
Let $\Pi:=\{(\lambda,\mu)~|~0\leq \lambda,\mu\leq p^n-2\}$. For each pair $(\lambda,\mu)\in \Pi$, we define a $k'\tilde{H}$-module $V_{\lambda,\mu}$: let $V_{\lambda,\mu}$ be a $1$-dimensional vector space over $k'$ on which $\left[ {\begin{array}{*{20}{c}}
  a&b \\
  0&d
\end{array}} \right]\in \tilde{H}$ acts as scalar multiplication by $a^{\lambda+\mu}d^\mu$. Then $\{V_{\lambda,\mu}~|~(\lambda,\mu)\in \Pi\}$ is a set of representatives of isomorphism classes of simple $k'\tilde{H}$-modules. Let $T_\lambda:={\rm Res}_{k'H}(V_{\lambda,0})$, then $\{T_\lambda~|~0\leq\lambda\leq p^n-2\}$ is a set of representatives of isomorphism classes of simple $k'H$-modules (see \S\ref{simplemodulesofSL}). So every simple $k'H$-module extends to a simple $k'\tilde{H}$-module.

Let $T_\lambda$ be a simple $k\tilde{H}c$-module. By \cite[Chapter 5, Lemma 5.8 (ii)]{NT} and \cite[Corollary 8.20]{Navbook}, there exist a simple $k'\tilde{H}c_1$-module $V_{\lambda_1,\mu_1}$, a simple $k'\tilde{H}c_t$-module $V_{\lambda_t, \mu_t}$, and a simple $k'(\tilde{H}/H)$-module $S$, such that
\begin{equation}\label{888}
{\rm Res}_{k'H}(V_{\lambda_1,\mu_1})=T_\lambda={\rm Res}_{k'H}(V_{\lambda_t,\mu_t})~~~{\rm and}~~~V_{\lambda_t,\mu_t}\cong V_{\lambda_1,\mu_1}\otimes_{k'}S. \tag{8.4.1}
\end{equation}
Note that we will also regard $S$ as a $k'\tilde{H}$-module or a $k'\Delta\tilde{H}$-module in an obvious way, and we still denote it by $S$. Then it is easy to see that
the functor sending a $k'\tilde{H}c_1$-module $V$ to the $k'\tilde{H}c_t$-module $V\otimes_{k'}S$ induces a Morita equivalence between $k'\tilde{H}c_1$ and $k'\tilde{H}c_t$. Assume that this Morita equivalence is induced by a $(k'\tilde{H}c_t,k'\tilde{H}c_1)$-bimodule $M'$. Then
$$M'\cong M'\otimes_{k'\tilde{H}c_1}k'\tilde{H}c_1\cong k'\tilde{H}c_1\otimes_{k'} S,$$
where the second isomorphism is by the definition of the functor above. So $M'$ is isomorphic to a direct summand of
$$k'\tilde{H}\otimes_{k'}S\cong {\rm Ind}_{\Delta\tilde{ H}}^{\tilde{H}\times\tilde{H}^{\rm op}}(S),$$
where the isomorphism is by \cite[Proposition 2.8.19]{Linckelmann}. Hence $M'$ has a vertex contained in $\Delta\tilde{H}$. Note that $S$ is a 1-dimensional $k'\tilde{H}$-module, hence the restriction of $S$ to any $p$-subgroup of $\Delta \tilde{H}$ is a trivial module. Using the Mackey formula, we easily see that $M'$ has trivial source. So $M'$ induces a splendid Morita equivalence.

If $S$ is defined over $k$, then $M'$ is defined over $k$, and then the proposition follows by \cite[Proposition 4.5 (c)]{Kessar_Linckelmann}. Our next task is to show that we can choose $S$ satisfying (\ref{888}), such that $S$ is defined over $k$.
We know that $c_1$ and $c_t$ are contained in $(k'C_{\tilde{H}}(P))^{\tilde{H}}$ (see \cite[Theorem 6.2.6]{Linckelmann}). We can easily calculate that $C_{\tilde{H}}(P)=\left\{\left[ {\begin{array}{*{20}{c}}
  a&b \\
  0&a
\end{array}} \right] \middle|~a\in \F_{p^n}^\times,~b\in \F_{p^n}\right\}= PZ(\tilde{H})$. Hence $k'C_{\tilde{H}}(P)$ has exactly $p^n-1$ blocks, and these blocks are exactly the blocks of $k'Z(\tilde{H})$.
Since $(k'C_{\tilde{H}}(P))^{\tilde{H}}\subseteq Z(k'C_{\tilde{H}}(P))$, $k'\tilde{H}$ has at most $p^n-1$ blocks. Denote by $\chi_{\lambda,\mu}:\tilde{H}\to {k'}^\times$ the character of $V_{\lambda,\mu}$. For an element
$u=\left[ {\begin{array}{*{20}{c}}
  a&0 \\
  0&a
\end{array}} \right]\in Z(\tilde{H})$, $\chi_{\lambda,\mu}(u)=a^{\lambda+2\mu}$. As $(\lambda,\mu)$ runs over $\Pi$, we see that the restriction $\chi_{\lambda,\mu}|_{Z(\tilde{H})}$ runs over all $p^n-1$ irreducible $k'$-characters of $Z(\tilde{H})$. It follows that $k'\tilde{H}$ has at least $p^n-1$ blocks (see e.g. \cite[Proposition 6.5.6]{Linckelmann} or \cite[Definition 3.1]{Navbook}). So the blocks of $k'\tilde{H}$ are exactly blocks of $k'Z(\tilde{H})$. Let $\xi\in \F_{p^n}\subseteq k'$ be a primitive $(p^n-1)$-th root of unity. Then $Z(\tilde{H})$ is a cyclic group generated by the element $z:=\left[ {\begin{array}{*{20}{c}}
  \xi&0 \\
  0&\xi
\end{array}} \right].$
Let $n_1:=\lambda_1+2\mu_1$ and $n_t:=\lambda_t+2\mu_t$. Then $\chi_{\lambda_1,\mu_1}(z)=\xi^{n_1}$, $\chi_{\lambda_t,\mu_t}(z)=\xi^{n_t}$.

By the discussion above, we know that $c_1$ is exactly a block of $k'Z(\tilde{H})$. Since $V_{\lambda_1,\mu_1}$ belongs to the block $c_1$ of $k'\tilde{H}$, the character $\chi_{\lambda_1,\mu_1}|_{Z(\tilde{H})}$ belongs to the block $c_1$ of $k'Z(\tilde{H})$. Hence we have
$$c_1=\frac{\chi_{\lambda_1,\mu_1}(1)}{|Z(\tilde{H})|}\sum_{i=0}^{p^n-2}\chi_{\lambda_1,\mu_1}(z^{-i})z^i
=\frac{1}{p^n-1}\sum_{i=0}^{p^n-2}\xi^{-in_1}z^i.$$
So $\F_p[c_1]=\F_p[\xi^{n_1}]$, and similarly we have $\F_p[c_t]=\F_p[\xi^{n_t}]$. Note that the order of the cyclic group generated by $\xi^{n_1}$ is $(p^n-1)/n_1$. Since $k=\F_p[c_1]=\F_p[c_t]$ (see Lemma \ref{filedvalue}), $\xi^{n_t}$ is equal to $(\xi^{n_1})^m$ for some positive integer $m<(p^n-1)/n_1$, such that $m$ is coprime to $(p^n-1)/n_1$.
So if $m$ is even, then $(p^n-1)/n_1$ is odd.

If $m$ is odd, we take $S=V_{0,\frac{(m-1)n_1}{2}}$. Note that $V_{\lambda_t,\mu_t}$ is uniquely determined by $V_{\lambda_1,\mu_1}$ and $S$.  By definition, $H$ acts trivially on $S$. Denote by $\chi_S:\tilde{H}\to {k'}^\times$ the character of $S$,
then $\chi_S(z)=\xi^{(m-1)n_1}$, $\chi_{\lambda_1,\mu_1}(z)\chi_S(z)=\xi^{n_t}$, and thus $V_{\lambda_1,\mu_1}\otimes_{k'}S$ is a simple $k'\tilde{H}c_t$-module. So $S$ satisfies (\ref{888}). Since ${\rm dim}_{k'}(S)=1$, $S$ is defined over $k$ if and only if the values of $\chi_S$ are contained in $k$.  It is easy to see that the values of $\chi_S$ are contained in the field $\F_p[\xi^{\frac{(m-1)n_1}{2}}]$, hence contained in $k=\F_p[\xi^{n_1}]$.

It is easy to check that $0\leq \frac{(m-1)n_1+(p^n-1)}{2}< p^n-1$. If $m$ is even, we take $S=V_{0,\frac{(m-1)n_1+(p^n-1)}{2}}$. By definition, $H$ acts trivially on $S$. Denote by $\chi_S:\tilde{H}\to {k'}^\times$ the character of $S$,
then $\chi_S(z)=\xi^{(m-1)n_1+(p^n-1)}=\xi^{(m-1)n_1}$, and $\chi_{\lambda_1,\mu_1}(z)\chi_S(z)=\xi^{n_t}$, and thus $V_{\lambda_1,\mu_1}\otimes_{k'}S$ is a simple $k'\tilde{H}c_t$-module. So $S$ satisfies (\ref{888}). It is easy to see that the values of $\chi_S$ are contained in the field $\F_p[\xi^{\frac{(m-1)n_1+(p^n-1)}{2}}]$. Note that $\frac{(m-1)n_1+(p^n-1)}{2}=\frac{(m-1+(p^n-1)/n_1)n_1}{2}$ and $m-1+(p^n-1)/n_1$ is even. So $\F_p[\xi^{\frac{(m-1)n_1+(p^n-1)}{2}}]$ is contained in $k=\F_p[\xi^{n_1}]$, and $S$ is defined over $k$. This completes the proof.   $\hfill\square$

\bigskip\noindent\textbf{Acknowledgements.} We are deeply indebted to a very kind, patient referee for a long list of very detailed comments and suggestions, which simplify the proofs of Proposition \ref{descent} and contribute to a significant improvement in the presentation of the paper. We wish to thank Joseph Chuang, Radha Kessar, Markus Linckelmann and Anderi Marcus for many very helpful conversations. The first author is grateful to Joseph Chuang for his instruction on the reference \cite{Oku} and for his proof of Proposition \ref{splendid}, and to Yuanyang Zhou for his guidance and constant support.

The second and third authors acknowledge the support by National Key R\&D Program of China (No. 2020YFE 0204200) and  NSFC (No. 11631001). The first author acknowledges former support by China Scholarship Council (No. 202006770016).



\end{document}